\documentclass[conference]{IEEEtran}

\usepackage{amsmath,amssymb,graphicx,algorithm,color,url}
\newcommand{\beqn}{\begin{equation}}
\newcommand{\eeqn}{\end{equation}}
\newcommand{\eref}[1]{(\ref{#1})}

\newtheorem{theorem}{Theorem}[section]

\newtheorem{defi}[theorem]{Definition}

\def \om {{\omega    }}
\def \N  {{\mathbb  N}}
\def \R  {{\mathbb  R}}
\def \cD {{\mathcal D}}
\def \cI {{\mathcal I}}
\def \cM {{\mathcal M}}
\def \cP {{\mathcal P}}
\def \cS {{\mathcal S}}

\begin{document}

\title{A Different View on the Vector-valued \\
       Empirical Mode Decomposition (VEMD)}

\author{\IEEEauthorblockN{Boqiang Huang}
\IEEEauthorblockA{Institut f\"ur Mathematik,
                  Universit\"at Paderborn, Germany \\
                  E-mail: bhuang@math.uni-paderborn.de}
\and
\IEEEauthorblockN{Angela Kunoth}
\IEEEauthorblockA{Mathematisches Institut,
                  Universit\"at zu K\"oln, Germany \\
                  E-mail: kunoth@math.uni-koeln.de}}

\maketitle

\begin{abstract}
  The empirical mode decomposition (EMD) has achieved its reputation
  by providing a multi-scale time-frequency representation of
  nonlinear and/or nonstationary signals. To extend this method
  to vector-valued signals (VvS) in multi-dimensional (multi-D)
  space, a multivariate EMD (MEMD) has been designed recently,
  which employs an ensemble projection to extract local extremum
  locations (LELs) of the given VvS with respect to different
  projection directions.  This idea successfully overcomes the
  problems of locally defining extrema of VvS. Different from the
  MEMD, where vector-valued envelopes (VvEs) are interpolated based on
  LELs extracted from the 1-D projected signal, the vector-valued EMD
  (VEMD) proposed in this paper employs a novel back projection
  method to interpolate the VvEs from 1-D envelopes in the projected
  space.  Considering typical 4-D coordinates (3-D location and
  time), we show by numerical simulations that the VEMD
  outperforms state-of-art methods.

\end{abstract}


\section{Introduction}

The empirical mode decomposition (EMD) was firstly designed for
nonlinear and/or nonstationary signal analysis \cite{Huang1998}.
Combined with the Hilbert transform, the Hilbert-Huang-Transform
(HHT) \cite{Huang2005} provides a finer time-frequency spectrum of a
given signal compared with other well-known methods such as
the Fourier transform, the wavelet transform, or the
  Wigner-Ville transform.  Moreover, the EMD does not require any
pre-defined basis. It decomposes a given signal $f(t)$ into a
  finite number of intrinsic mode functions (IMFs), $f_j(t):=a_j(t)
\cos(\theta_j(t)), j=1,\ldots,J$, and a monotonic trend $r_{J+1}(t)$,
i.e.,
\beqn\label{EMD}
  f(t) := \sum_{j=1}^J f_j(t) +  r_{J+1}(t),\qquad t\in [0,T].
\eeqn
Here, each $f_j(t)$ might be considered as an amplitude-modulated and
frequency-modulated signal, or as a generalized Fourier component
\cite{Hou2011}. The properties of the signal model are extensively
studied in \cite{Huang1998, Huang2013} and the references therein.
Based on \eref{EMD}, the instantaneous frequency of each IMF is well
defined by $\om_j(t):= \theta'_j(t)$ \cite{Huang2009}. In
  addition, the corresponding Hilbert amplitude spectrum can be
naturally derived as $\{a_j(t) \mbox{ on the curves } (t,
\omega_j(t)), ~ t\in [0,T], ~j=1,\ldots,J\}$ \cite{Huang2013}. Up to
now, the EMD and its many variations have been successfully
employed in many disciplines, such as signal processing
\cite{Huang2005}, hydrology \cite{Rudi2010}, and geophysics
\cite{Huang2008}.

In general, the EMD aims at sequentially extracting each IMF
$f_j(t)$ through a filtering operation called sifting process.
Considering the signal at the $k$th level,
$x_k(t) := f(t)-\sum_{j=1}^{k-1} f_j(t)$, the sifting operator
$\cS$ can recursively be defined by
$\cS^n [x_k](t) := \cS^{n-1} [x_k](t) - \cM[\cS^{n-1} [x_k]](t)$,
where $n$ is the iteration number, $\cS^0 [x_k](t) := x_k(t)$,
and $\cM [\cdot]$ represents the local trend approximation
operator. In the classic EMD \cite{Huang1998}, $\cM [x](t)$ is
defined as the mean curve of the upper and lower envelopes
which are defined by  cubic spline interpolation of
the local maxima and local minima of $x(t)$, respectively. The sifting
process stops at $n = N$ until some \textit{ad hoc} criterion
is met, e.g. $f_k(t) := \cS^N [x_k](t)$ mimics in some sense a
generalized Fourier component.

The EMD is a completely data-driven decomposition method which
heavily depends on the definition of local extrema. Since their
  definition is unclear in higher dimensions, or the extrema are
  nonunique, it is, therefore, difficult to extend the method to these
  cases.  The original signal or univariate time series decomposed in
  (\ref{EMD}) consisted of data of the form $f:[0,T]\to \R$, mapping
one-dimensional data in one-dimensional space. In higher dimensions,
we distinguish between {\em multivariate} and {\em vector-valued}
data. The {\em multivariate case} consists of `cube' data, e.g., an
image in 2-D, or a volume in 3-D. In this case, the function is of the
form $f:[0,T]^d\to \R$ with $d=2,3$ being the space dimension. Here,
a local extremum may be defined as the strict extremum in a
pre-defined neighborhood \cite{Jager2010}. However, in the {\em
    vector-valued} case, the data are multi-D `curve' data which are
  of the form $f:[0,T]\to \R^d$. In this case, the definition of local
  extrema is much more complicated since even the notion of
  `neighborhood' is unclear.

In this paper, we concentrate on the latter, more difficult, case. We
define the vector-valued signal as follows.
\begin{defi}\label{DefVSignal}
Given a finite number of 1-D signals
$f^{[i]}(t) \in C(\R) \cap L_{\infty}(\R), i=1,\ldots,d,~ t\in [0,T]$,
the corresponding vector-valued signal (VvS) is defined as
$F(t):=(f^{[1]}(t), \cdots, f^{[d]}(t))^T$, i.e., $F:[0,T]\to \R^d$.
\end{defi}

To develope an EMD-like decomposition for VvS, a straightforward
idea is to decompose a complex signal, $F(t):=(f^{[1]}(t),
f^{[2]}(t))^T$, by applying the classic EMD to the real and imaginary
parts separately. However, this usually leads to a different
number of IMFs for each of the two components which is an
  undesired effect \cite{Tanaka2007}.  In fact, if we reconsider the
signal decomposition model in 3-D (2-D location and time), and then
define the corresponding $j$th complex IMF as $F_j(t):=(f^{[1]}_j(t),
f^{[2]}_j(t))^T = (a_j^{[1]}(t) \cos(\theta_j^{[1]}(t)),
a_j^{[2]}(t)\cos(\theta_j^{[2]}(t))^T$, we will observe an interesting
fact that each $F_j(t)$ is a rotation invariant component. Its
3-D envelope $(a_j^{[1]}(t), a_j^{[2]}(t), t)$ should be some
tube-shaped surface tightly enclosing the 3-D curve $F_j(t)$.
Then its 3-D local trend/mean $\cM[F_j](t)$ can be considered as the
barycenter curve of the 3-D envelope surface.  This observation
implies a possible way to approximate $\cM[F_j](t)$ from the view of
statistics. This means that we should interpolate the envelopes
based on the local extrema of $F_j(t)$ along a selected projection
direction, and then average all interpolated envelopes with
respect to all possible projection directions, following a idea
  from \cite{Rilling2007}.  Here, we define within a general concept
the multi-D local extremum as follows.

\begin{defi}\label{DefMultiDExtremum}
Given a VvS $F(t)$ and a selected unit projection direction $p$,
the local extremum of $F(t)$ along $p$ is defined as the hyperpoint
$(F(t_k), t_k)$, where $t_k$ is the corresponding local extremum
location (LEL) of the projected 1-D signal $\cP_p[F](t)$.
\end{defi}

If the projection number approaches infinity, $\cM[F_j](t)$
should be the expectation of the local mean approximation based
on the projection. Such operations belong to a typical ensemble
approach. This implies that the uniform direction sampling
scheme should be an optimal choice \cite{Rehman2010, Rehman2010_2}.
Moreover, decomposition tests with white noise show that such
EMD-like decomposition of VvS perfectly inherit the dyadic filter
bank property of the classic EMD \cite{Rehman2011, Mandic2013}.

To further study the above mentioned local mean approximation, we
consider the vector-valued envelope (VvE) interpolation in a
different but more general way than in
\cite{Rilling2007, Rehman2010, Rehman2010_2}. Here, the methods directly
interpolate the VvE in the multi-D space and assume that the following
property is satisfied: if we project the interpolated VvE
along the corresponding projection direction, the projected curve
should be nothing but the envelope interpolated based on the local
extrema in the projected 1-D space. Unfortunately, this assumption
only meets the requirement of the direct interpolation method
(without constraints), e.g., the cubic spline interpolation (CSI),
but might be defective for the others, e.g. optimization based
method \cite{Huang2013}, when additional constraints have to be
maintained. Thus, a more natural way to obtain the VvE is to
firstly interpolate the 1-D envelope in the projected space with the
desired method, and then back-project it into the original multi-D
space.

In this paper, we explain how to realize the
vector-valued EMD (VEMD) with an optimization based back projection.
To simplify the discussion, a typical signal decomposition model in
4-D space (3-D location and time) will be studied. In the following,
Section II recalls the recently developed multivariate EMD (MEMD) with
approximately uniform
direction sampling, Section III introduces the proposed VEMD,
Section IV presents the numerical studies, and Section V concludes
the paper.

\section{Multivariate empirical mode decomposition}

To keep the terminology from previous papers, in this section, the
MEMD concerns the VvS decomposition problem. For a concrete 4-D
decomposition problem, the signal model \eref{EMD} can be generalized
as
\beqn\label{VEMD}
    \begin{aligned}
        F(t) &:= \sum_{j=1}^J F_j(t) +  R_{J+1}(t),  \\
        \left( \begin{array}{c}
                 f^{[1]}(t)  \\
                 f^{[2]}(t)  \\
                 f^{[3]}(t)  \\
        \end{array}  \right) &:= \sum_{j=1}^J
        \left( \begin{array}{c}
                 f^{[1]}_j(t)  \\
                 f^{[2]}_j(t)  \\
                 f^{[3]}_j(t)  \\
        \end{array}  \right) +
        \left( \begin{array}{c}
                 r^{[1]}_{J+1}(t)  \\
                 r^{[2]}_{J+1}(t)  \\
                 f^{[3]}_{J+1}(t)  \\
        \end{array}  \right),
    \end{aligned}
\eeqn
where $F(t) := (f^{[1]}(t), f^{[2]}(t), f^{[3]}(t))^T$ is the given
4-D curve, $F_j(t) := (f^{[1]}_j(t), f^{[2]}_j(t), f^{[3]}_j(t))^T$
is the $j$th decomposed IMF, and
$R_{J+1}(t) := (r^{[1]}_{J+1}(t), r^{[2]}_{J+1}(t), r^{[3]}_{J+1}(t))^T$
is the monotonic trend.

\subsection{Approximately uniform sampling on a unit sphere}

To obtain equidistributed direction on a unit sphere, in this paper,
we employ a low-discrepancy sampling scheme based on transformed
Hammersley points, which can produce more uniform samples on the
hypersphere than other methods, like the polar coordinate lattices
or rotation method \cite{Cui1997}.

\begin{figure}[t]
\begin{minipage}[b]{.48\linewidth}
  \centering
  \centerline{\includegraphics[width=4.4cm]{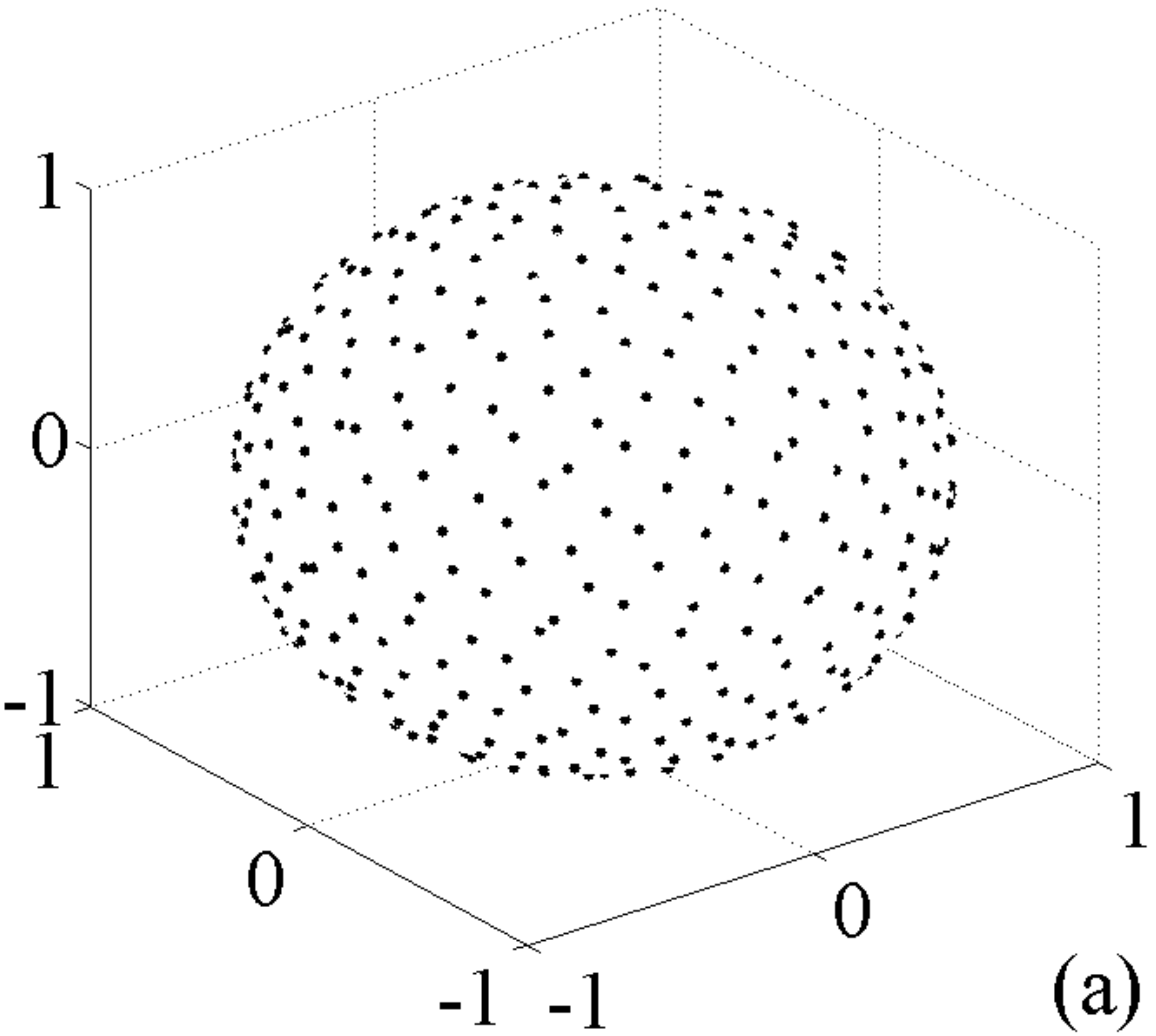}}
\end{minipage}
\hfill
\begin{minipage}[b]{.48\linewidth}
  \centering
  \centerline{\includegraphics[width=4.4cm]{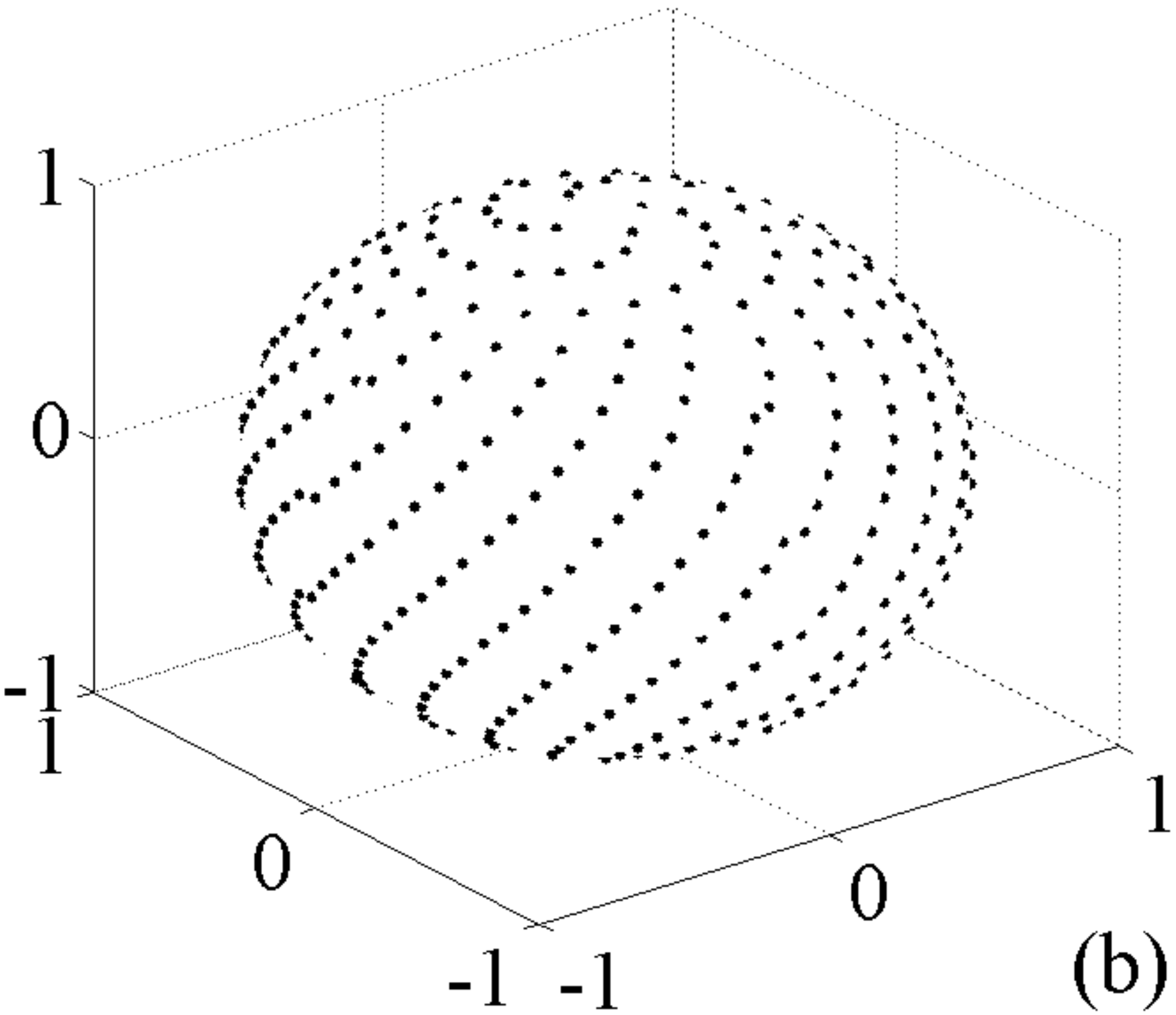}}
\end{minipage}
\caption{\small Two sampling examples based on Hammersley points.
(a) $p_m \in \cI_2^{512}$, (b) $p_m \in \cI_5^{512}$.}
\label{fig:Hammersley}
\end{figure}

Simply speaking, for any assumed direction number $M$ and any prime
base $b, b \geq 2, b \in \N$, a selected direction index
$m, m \in \{0,1,\cdots, M-1\},$ has a unique digit expansion
\beqn\label{DigitExpan}
  m = \sum_{j=0}^{r_m} c_{m_j} b^j,
\eeqn
where each $c_{m_j}$ is an integer in $[0, b-1]$. The Van der Corput
sequence can be defined as
\begin{small}
\beqn\label{VanDerCorput}
  \{ z_b(m) | z_b(m) := \sum_{j=0}^{r_m} c_{m_j} b^{-j-1}, m = 0,1,\cdots, M-1 \},
\eeqn
\end{small}
and the corresponding Hammersley points set is \cite{Niederreiter1992}
\beqn\label{Hammersley}
 \cI_b^M := \left\{ \left( \frac{m}{M}, z_b(m) \right), m = 0,1,\cdots, M-1 \right\}.
\eeqn
Since $\frac{m}{M} \in [0, 1)$, and $z_b(m) \in [0, 1)$, we can
transform the equidistributed square $[0, 1) \times [0, 1)$ into
equidistributed cylinder $[0, 2\pi) \times [-1, 1)$ by setting
$\phi := 2 \pi \frac{m}{M}$ and $z := 2 z_b(m) - 1$. Finally the
uniformly sampled projection directions are
$p_m := (\sqrt{1-z^2} \cos(\phi), \sqrt{1-z^2}) \sin(\phi), z)^T, m=0,1,\cdots, M.$
Fig.~\ref{fig:Hammersley} shows two sampling examples for $M = 512$,
and $b=2, 5$ respectively.

It should be noted that a more satisfactory sampling scheme in
multi-D space can be generated by introducing different prime bases
which connects to the Hammersley points based on Halton sequences
\cite{Niederreiter1992}.

\subsection{Multivariate empirical mode decomposition}\label{SecMEMD}

The recent MEMD method has been shown its strength for VvS decomposition,
especially for mode alignment problem and noise-assisted applications.
The white noise decomposition test illustrates its remarkable dyadic
filter bank property comparing to classic EMD and ensemble EMD
\cite{Mandic2013}. The whole algorithm can be summarized in
Algorithm~\ref{Alg_MEMD}.

\begin{algorithm}
\renewcommand{\thealgorithm}{1}
\caption{: MEMD in 4-D space}
\small
\label{Alg_MEMD}
\setlength{\IEEEilabelindent}{3pt}
\setlength{\IEEEiednormlabelsep}{2pt}
\setlength{\IEEEiedtopsep}{0pt}
\begin{IEEEitemize}
\item[1:] set $R_1(t) = F(t)$, $j = 1$, and the values of $M$ and $b$;
          generate the projection direction set $\cI_b^M$;
\item[2:] extract $j$th IMF by sifting process $(F_j(t) := \cS^N[R_j](t))$
\begin{IEEEitemize}
\item[(A)] set intermediate curve $H_{j,0}(t) = R_j(t)$ and $n = 0$;
\item[(B)] approximate the local mean $T_{j,n}(t) := \cM[H_{j,n}](t)$;
\begin{IEEEitemize}
\item[(a)] for every direction $p_m \in \cI_b^M, m=0,1,\cdots,M-1$, project
           $H_{j,n}(t)$ into 1-D space,
           $h_{j,n}^{p_m}(t):= \cP_{p_m}[H_{j,n}](t)$;
\item[(b)] detect all LELs of $h_{j,n}^{p_m}(t)$ (maximum location $t_{j,n,k^+}^{p_m}$
           and minimum location $t_{j,n,k^-}^{p_m}$, $k^+, k^- \in \N$) ;
\item[(c)] generate upper $U_{j,n}^{p_m}(t)$ and lower $V_{j,n}^{p_m}(t)$
           VvEs by interpolating data pairs
           $(H_{j,n}(t_{j,n,k^+}^{p_m}), t_{j,n,k^+}^{p_m})$ and
           $(H_{j,n}(t_{j,n,k^-}^{p_m}), t_{j,n,k^-}^{p_m})$ based on CSI;
\item[(d)] set $T_{j,n}(t) := \frac{1}{2M}\sum_{m=0}^{M-1} U_{j,n}^{p_m}(t) + V_{j,n}^{p_m}(t)$;
\end{IEEEitemize}
\item[(C)] update $H_{j,n+1}(t) = H_{j,n}(t) - T_{j,n}(t)$ and $n = n + 1$;
\item[(D)] calculate stopping criterion $\text{SD}_{j,n}$ as the one in \cite{Rilling2003};
\item[(E)] repeat steps (B) to (E) until $\text{SD}_{j,n} \leq \text{SD}_{\text{Thr}}$,
           or the max iteration number $N$ is met; define
           $j$th IMF $F_j(t) := H_{j,n}(t)$;
\end{IEEEitemize}
\item[3:] update $R_{j+1}(t) = R_j(t) - F_j(t)$, and $j = j + 1$;
\item[4:] repeat steps 2 to 3 until number of extrema in $R_j(t)$ is less than 2 or
          an expected IMF index is met, i.e., $j = J$.
\end{IEEEitemize}
\end{algorithm}

In MEMD, the local mean approximation $\cM[\cdot](t)$ could be
considered as an ensemble approach, which means the idea mean curve
$\hat{T}(t)$ may be exactly approximated as $M \rightarrow \infty$.
In other words, we could re-define the local mean by
$\hat{T}(t) := \cM[H](t):= \text{E}(\frac{1}{2}(U^{p_m}(t) + V^{p_m}(t))$,
where E denotes the expectation over the direction $p_m$. In real
application, the projection number $M$ needs not to be a large number
compromising both computational complexity and approximation performance.

\section{Vector-valued empirical mode decomposition} \label{Sec:VEMD}

In this section, a different EMD extension, namely the VEMD, will be
explained, in which the VvEs of the given VvS are generated by
back-projecting the 1-D envelopes interpolated in the projected space
to the original 4-D space. Since the naive back-projection might result
in infinite solutions, a novel optimization scheme will be designed to
guarantee an unique VvE which should be also as smooth as possible.

\subsection{Optimization based back projection} \label{Sec:VEMD_BP}

Keeping the notation system in algorithm \ref{Alg_MEMD} but ignoring
the trivial subindexes, in this subsection, $H(t), p_m$ denote the
considering VvS and the projection direction, $U_{p_m}(t), V_{p_m}(t)$
denote the upper and lower VvEs w.r.t $p_m$ in the original 4-D space,
and $u_{p_m}(t), v_{p_m}(t)$ denote the upper and lower envelopes
w.r.t $p_m$ in the 1-D projected space.

Assuming $H(t)$ and $p_m$ are given, and the $\cP_{p_m}[H](t)$ is the
projected 1-D signal, the local maxima/minima of $\cP_{p_m}[H](t)$ can
be easily detected at the corresponding time locations $t_{k^+}$ and
$t_{k^-}$, $k^+, k^- \in \N$. With some interpolation method, e.g. the
CSI, the 1-D upper $u_{p_m}(t)$ and lower $v_{p_m}(t)$ envelopes can be
interpolated based on the pairs $(\cP_{p_m}[H](t_{k^+}), t_{k^+})$ and
$(\cP_{p_m}[H](t_{k^-}), t_{k^-})$ separately. Now, we aim to project
$u_{p_m}(t)$ and $v_{p_m}(t)$ back to the original 4-D space in order to
obtain the VvEs
$U_{p_m}(t):=(u_{p_m}^{[1]}(t), u_{p_m}^{[2]}(t), u_{p_m}^{[3]}(t))^T$
and $V(t):=(v_{p_m}^{[1]}(t), v_{p_m}^{[2]}(t), v_{p_m}^{[3]}(t))^T$
w.r.t the direction $p_m$.

Mathematically speaking, such back projection problem is equivalent to
the solution problem of the linear system
$\cP_{p_m}[U_{p_m}](t)=u_{p_m}(t)$
(or $\cP_{p_m}[V_{p_m}](t)=v_{p_m}(t)$),
which has infinite solutions because the number of unknowns is larger
than the number of equations. However, the infinite solutions can be
constrained to unique one if next two properties can be maintained
simultaneously: a) the back-projected VvEs $U_{p_m}(t)$ (or
$V_{p_m}(t)$) should pass through the LELs $(H(t_{k^+}), t_{k^+})$
and $(H(t_{k^-}), t_{k^-})$ w.r.t direction $p_m$; b) the VvEs
$U_{p_m}(t)$ (or $V_{p_m}(t)$) should be as smooth as possible. To
investigate the smoothness of an unknown 4-D VvS
$X(t):=(x^{[1]}(t), x^{[2]}(t), x^{[3]}(t))$, we may employ an $n$th
order Sobolev norm funtional, i.e., the $L_2$ norm of the $n$th (weak)
derivative of the function $X$
\beqn\label{Norm}
  \resizebox{.9\hsize}{!}
  {$\boldsymbol{S}^{(n)}[X]
    := \| \cD^{(n)}[x^{[1]}] \|_{L_2}^2 + \| \cD^{(n)}[x^{[2]}] \|_{L_2}^2
                                         + \| \cD^{(n)}[x^{[3]}] \|_{L_2}^2 $},
\eeqn
where $\cD^{(n)}$ is the $n$th order derivative operator as its matrix
form for discrete problem is well-known .

Now, the back-projected VvEs can be uniquely obtained by solving
the following optimization problems
\begin{eqnarray*}
(P1) \quad \mbox{Minimize} &\;& \boldsymbol{S}^{(n)}[U_{p_m}](t) \\
         \mbox{~\,subject to} &\;& \cP_{p_m}[U_{p_m}](t)=u_{p_m}(t) \\
                           &\;& U_{p_m}(t_{k^+})=H(t_{k^+}).
\end{eqnarray*}
\begin{eqnarray*}
(P2) \quad \mbox{Minimize} &\;& \boldsymbol{S}^{(n)}[V_{p_m}](t)  \\
         \mbox{~\,subject to} &\;& \cP_{p_m}[V_{p_m}](t)=v_{p_m}(t)  \\
                           &\;& V_{p_m}(t_{k^-})=H(t_{k^-})
\end{eqnarray*}
In fact, (P1) and (P2) are both quadratic optimization problems with equality
constraints. They can be written as quadratic optimization problems
without constraints by solving the linear constraint system and then
implementing the variable reduction.

\subsection{Vector-valued empirical mode decomposition}

Comparing to the MEMD method in section \ref{SecMEMD}, the VEMD employs
the back projection to obtain the VvEs from the envelopes interpolated
in the 1-D projected space w.r.t a selected projection direction.

Most of the computation steps in VEMD algorithm are as the same as the
ones in algorithm \ref{Alg_MEMD} except the step 2:(B):(c) which should
be replaced by

\begin{algorithm}
\renewcommand{\thealgorithm}{2}
\caption{: VEMD in 4-D space}
\small
\label{Alg_VEMD}
\setlength{\IEEEilabelindent}{3pt}
\setlength{\IEEEiednormlabelsep}{2pt}
\setlength{\IEEEiedtopsep}{0pt}
\begin{IEEEitemize}
\item[2:] (B):(c) generate upper $u_{j,n}^{p_m}$ and lower $v_{j,n}^{p_m}(t)$
  envelope in the 1-D projected space by interpolating data pairs
  $(\cP_{p_m}[H_{j,n}](t_{j,n,k^+}^{p_m}), t_{j,n,k^+}^{p_m})$ and
  $(\cP_{p_m}[H_{j,n}](t_{j,n,k^-}^{p_m}), t_{j,n,k^-}^{p_m})$ based on
  CSI; then solve the corresponding back projection optimization problem
  (P1) and (P2) to obtain the 4-D upper $U_{j,n}^{p_m}(t)$ and lower
  $V_{j,n}^{p_m}(t)$ envelopes;
\end{IEEEitemize}
\end{algorithm}

It should be noted that in algorithm \ref{Alg_VEMD}, the envelope
interpolation in the 1-D projected space can be implemented by any
reasonable interpolation method in order to meet different mathematical
requirements, e.g. the method in \cite{Huang2013}.

\section{Numerical studies}

In this section, all the simulations are implemented in MATLAB on a
laptop equipped with an i7-4700 quad-core CPU, 8 GB memory and under
Windows. All optimization problems
are solved by the standard CVX toolbox from \cite{Grant2011}. To distinguish
the associated variables in MEMD and VEMD, we employ the subscripts M and V.

In the VEMD, the order of the derivative operator in \eref{Norm}, the
projection direction number and the prim base in \eref{Hammersley}
are free parameters, each of which may effect the behavior of the
method. To determine each one, we study the decomposition problem of
the following VvS
\begin{small}
\beqn\label{Simu1}
    \begin{aligned}
        F(t) &:= X(t) + Y(t), \; t \in [0, 1],\\
        X(t) &:=
        \left( \begin{array}{c}
            x^{[1]}(t)  \\
            x^{[2]}(t)  \\
            x^{[3]}(t)  \\
        \end{array}  \right) =
        \left( \begin{array}{c}
            (1 + \cos(2 \pi t)) \sin(20 \pi t)  \\
            (2 + \cos(4 \pi t)) \sin(20 \pi t)  \\
            (3 + \cos(6 \pi t)) \sin(20 \pi t)  \\
        \end{array}  \right), \\
        Y(t) &:=
        \left( \begin{array}{c}
            y^{[1]}(t)  \\
            y^{[2]}(t)  \\
            y^{[3]}(t)  \\
        \end{array}  \right) =
        \left( \begin{array}{c}
              \sin(4 \pi t)  \\
            2 \sin(4 \pi t)  \\
            3 \sin(4 \pi t)  \\
        \end{array}  \right).
    \end{aligned}
\eeqn
\end{small}
Comparing to the signal model in \eref{VEMD}, the components $X(t)$
and $Y(t)$ can be considered as the first IMF $F_1(t)$ and the residual
$R_2(t)$, respectively. The reason that we set a common frequency of each
$x^{[i]}(t)$ and $y^{[i]}(t)$ is to facilitate the discussion. Since the
MEMD, as the same as the VEMD, has the mode alignment property (see
\cite{Mandic2013}), different frequencies involved in the components $x^{[i]}(t)$
or $y^{[i]}(t)$ undoubtedly lead to many more decomposed IMFs, each of which
should contain a particular frequency.

\subsection{Derivative order}

\begin{figure}[t]
\centering
\includegraphics[width=8.8cm]{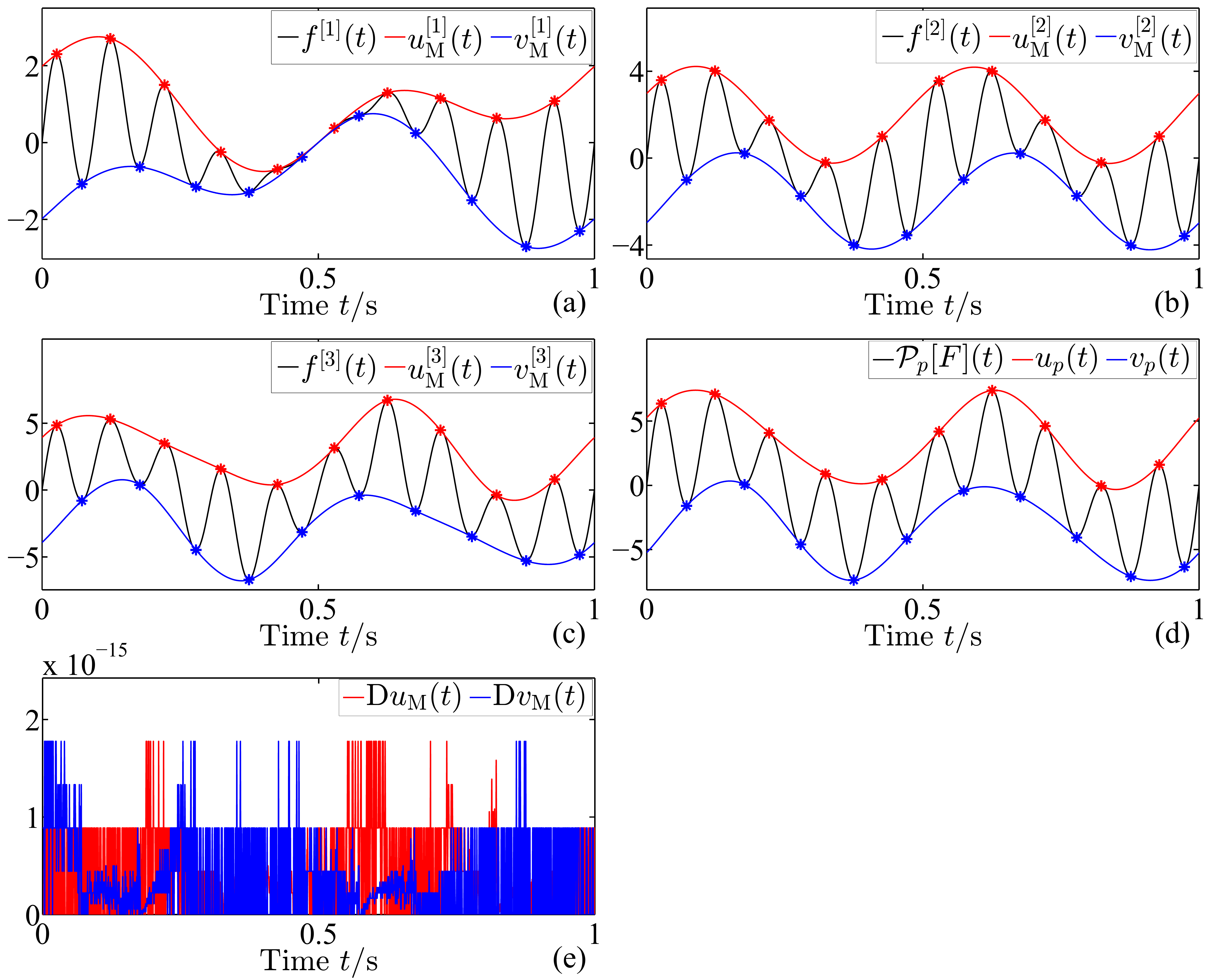}
\caption{\small
The envelope interpolation in MEMD. (a)-(c) The VvS $F(t)$ and its
interpolated VvEs $U_{\text{M}}(t)$ and $V_{\text{M}}(t)$ w.r.t an unit
projection direction $p$; (d) The projected signal $\cP_p[F](t)$ and
its interpolated envelopes $u_p(t)$ and $v_p(t)$; (e) The absolute
difference between the projected VvEs and the envelopes in the projected
space, e.g. $\text{D}u_{\text{M}}(t):=|\cP_p[U_{\text{M}}](t)-u_p(t)|$.}
\label{fig:MEMDProj}
\end{figure}

Selecting a unit projection direction
$p=[\frac{1}{2}, \frac{1}{2}, \frac{\sqrt{2}}{2}]$, we can easily obtain
the projected VvS $\cP_p[F](t)$ together with its LELs $t_{k^+}/t_{k^-}$
(Fig.\ref{fig:MEMDProj} (d)). Based on these locations, by using a CSI,
we can interpolate the envelopes $u_p(t), v_p(t)$ in the projected space
(Fig.\ref{fig:MEMDProj} (d)), and the VvEs
$U_{\text{M}}(t):=(u_{\text{M}}^{[1]}(t),u_{\text{M}}^{[2]}(t),u_{\text{M}}^{[3]}(t))^T$
and
$V_{\text{M}}(t):=(v_{\text{M}}^{[1]}(t),v_{\text{M}}^{[2]}(t),v_{\text{M}}^{[3]}(t))^T$
in the original 4-D space (Fig.\ref{fig:MEMDProj} (a)-(c)).
Fig.\ref{fig:MEMDProj} (e) illustrates that the projected VvEs
$(\cP[U_M](t),\cP[V_M](t))$ are nothing but the envelopes $(u_p(t),v_p(t))$
in the projected space. In other words, Fig.\ref{fig:MEMDProj} graphically
explains why the MEMD interpolates the VvEs directly in the original space
but not in the projected space.

\begin{figure}[t]
\centering
\includegraphics[width=8.8cm]{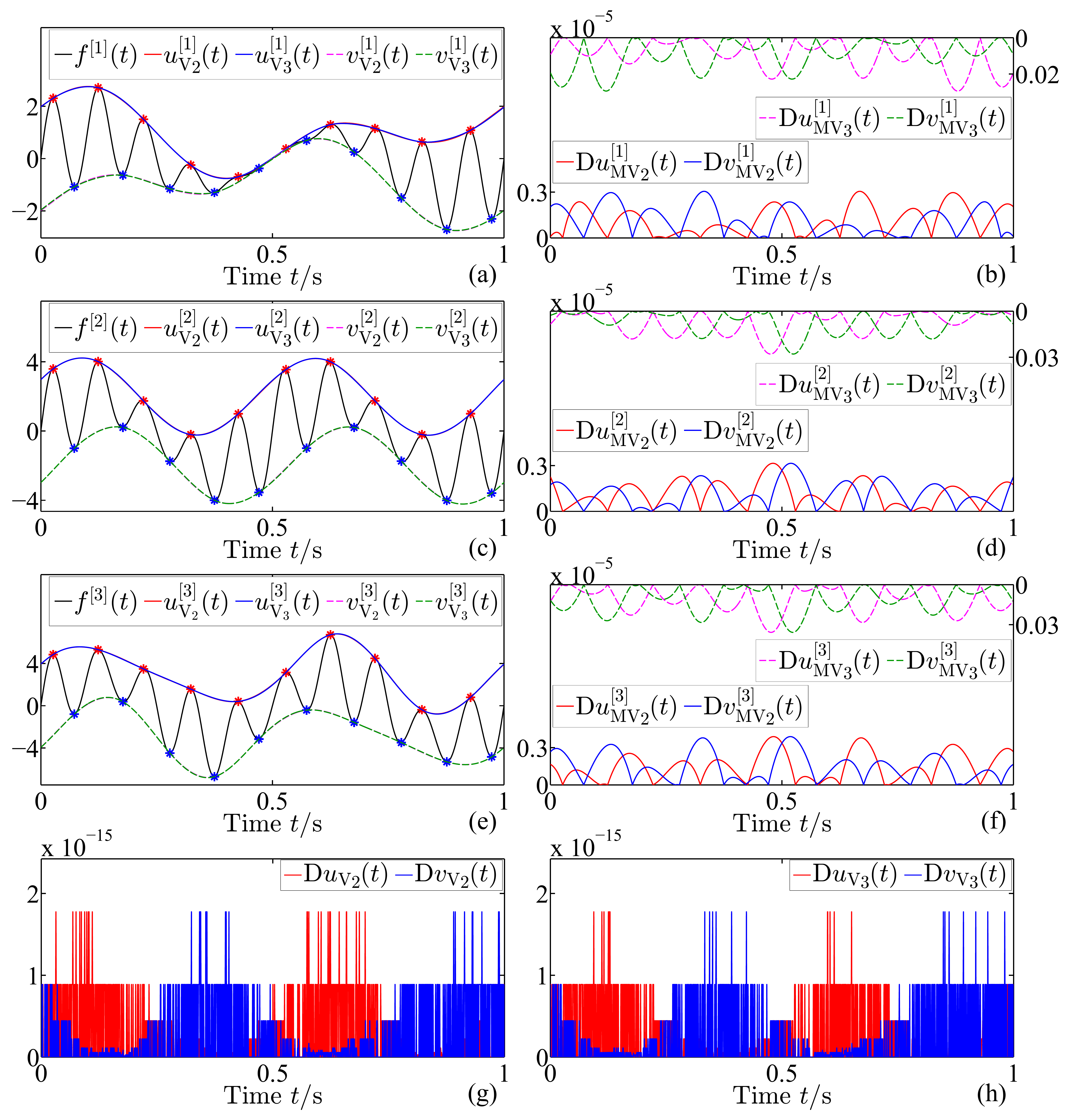}
\caption{\small
The envelope interpolation in VEMD. Left column (a)(c)(e): the VvS and
its VvEs obtained by solving back projection problems (P1) and (P2)
w.r.t an unit projection direction $p$. Subscripts $\text{V}_2$ and
$\text{V}_3$ of the solution denote the 2nd and 3rd order derivative
selected in \eref{Norm} for VvE interpolation; Right column (b)(d)(f):
absolute difference between interpolated VvEs in VEMD and the ones in
MEMD shown in Fig.\ref{fig:MEMDProj}, e.g.
$\text{D}U_{\text{MV}_2}(t):=|U_{\text{M}}(t)-U_{\text{V}_2}(t)|$; Last
row (g)(h): The absolute difference between the projected VvEs and the
envelope in the projected space, e.g.
$\text{D}u_{\text{V}_2}(t):=|\cP_p[U_{\text{V}_2}](t)-u_p(t)|$.}
\label{fig:VEMDDerivN}
\end{figure}

Now, based on the interpolated envelopes together with the LELs in the
projected space, we can obtain the 4-D VvEs using the proposed back
projection method, i.e. solving the optimization problems (P1) and (P2)
described in Sec.\ref{Sec:VEMD_BP}. Considering that, normally, the
first order derivative won't be an optimum choice for smoothness
measurement \cite{Huang2013}, we set the 2nd or 3rd order derivative
in \eref{Norm} alternatively for simulation. Fig.\ref{fig:VEMDDerivN}
presents the corresponding solved VvEs based on the data shown in
Fig.\ref{fig:MEMDProj} (d).

In Fig.\ref{fig:VEMDDerivN}, sub-figures (a)(c)(e) illustrate
back-projected VvEs $U_{\text{V}_2}(t), V_{\text{V}_2}(t)$ and
$U_{\text{V}_3}(t), V_{\text{V}_3}(t)$, where the subscripts
$\text{V}_2$ and $\text{V}_3$ denote the 2nd and 3rd order derivative
selected in VEMD each. To evaluate the interpolation performance, we
take the VvEs interpolated in MEMD as a benchmark.
Sub-figures (b)(d)(f) present the absolute differences between each
VvE in VEMD and the corresponding VvE in MEMD, e.g.
$\text{D}U_{\text{MV}_2}(t) := (\text{D}u_{\text{MV}_2}^{[1]}(t),
 \text{D}U_{\text{MV}_2}^{[2]}(t), \text{D}U_{\text{MV}_2}^{[3]}(t))^T
 := |U_{\text{M}}(t)-U_{\text{V}_2}(t)|$.
These figures imply that the VEMD with 2nd order derivative would
be much close to the MEMD comparing to the VEMD with 3rd order
derivative. This interesting phenomenon can be effortlessly understood
if we can recall that the CSI requires the interpolated curve to be at
most 2nd order continuously differentiable. Therefore, the VEMD with
2nd order derivative can produce similar results to the ones from MEMD,
while the VEMD with 3rd order derivative can provide smoother
interpolated curve. Finally, sub-figures (g)(h) imply both VvEs
interpolated in VEMD satisfy the first equality constraint in (P1)
and (P2) with high accuracy.

\subsection{Projection number and prime base}

In section \ref{SecMEMD} we have shown that, given a fixed prime
base for projection sampling, the local mean approximation operator
$\cM[\cdot]$ might be considered as an ensemble approach. In practice,
we have to determine a finite projection number such that the
approximation would be good enough or at least approximately
convergent. Let's consider again the VvS $F(t)$ in \eref{Simu1}. If
we apply the operator on $F(t)$ with a fixed projection number $M$,
$\cM^M[F](t)$ shall be the approximated local mean, e.g.
$\cM_{\text{M}}^M[F](t)$, the result of step (B) in the MEMD
(algorithm \ref{Alg_MEMD}), or the $\cM_{\text{V}_2}^M[F](t)$ or
$\cM_{\text{V}_3}^M[F](t)$ in the VEMD (algorithm \ref{Alg_VEMD}).
To measure the approximation efficiency, we define a percent root
mean squared difference (PRD) as an error function of the projection
number $M$,
\beqn\label{PRD}
  \text{PRD}(M) := \frac{\|\cM^M[F](t)-Y(t)\|_{L_2}}{\|Y(t)\|_{L_2}} \times 100\%.
\eeqn

\begin{figure}[t]
\centering
\includegraphics[width=8.8cm]{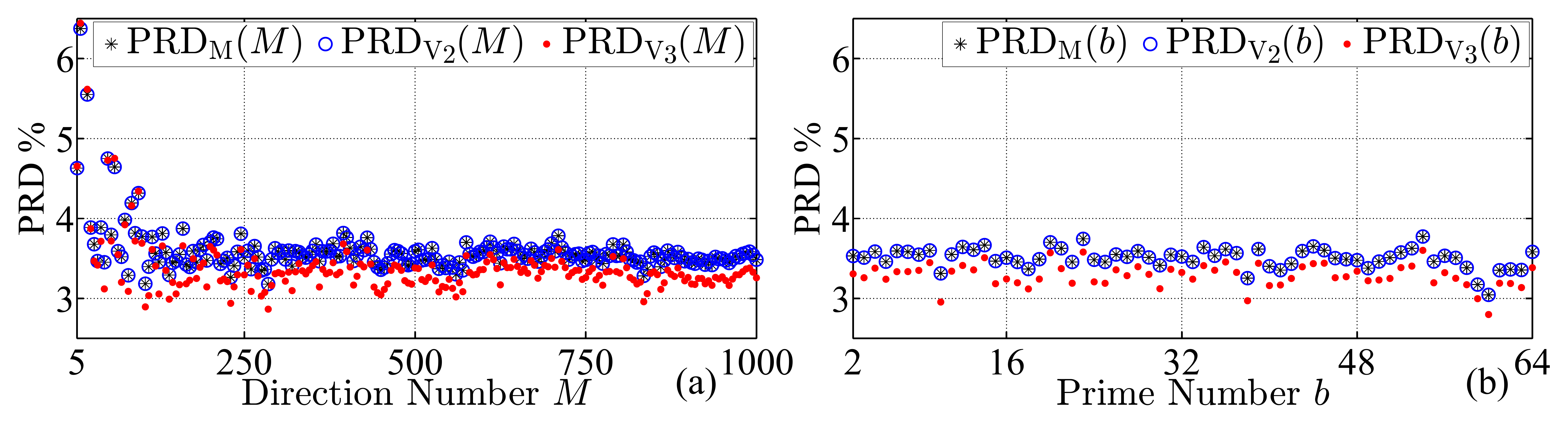}
\caption{\small
Local mean approximation performance of operators $\cM_{\text{M}}^M[\cdot]$,
$\cM_{\text{V}_2}^M[\cdot]$ and $\cM_{\text{V}_3}^M[\cdot]$ in MEMD and VEMD.
(a) PRD values of the three operators for $b = 2$, and $M \in [5, 1000]$;
(b) PRD values of the three operators for $M = 512$, and $b \in [2, 80]$.}
\label{fig:DirectPrime}
\end{figure}

Fig.\ref{fig:DirectPrime} (a) describes the approximation
performance of the operators $\cM_{\text{M}}^M[\cdot]$,
$\cM_{\text{V}_2}^M[\cdot]$ and $\cM_{\text{V}_3}^M[\cdot]$ for
$b = 2, M \in [5, 1000]$. As can be seen, all PRD sequences decay
dramatically in the beginning, and each one might be approximately
convergent when $M \geq 256$. Since the sampled projection directions
are not ideally uniform distributed on the sphere \cite{Cui1997},
the sampling error may result in slight oscillations on the PRD
sequence. Therefore, we set $M = 512$ for rest simulations.

Fig.\ref{fig:DirectPrime} (b) presents the performance of all
three approximation operators with different prime bases
$b \in [2, 80]$. It shows that different prime numbers do not
effect the PRD values significantly. This is because the prime
base only effects on the projection locations on the sphere
(see Fig.\ref{fig:Hammersley}). When the projection number is
large enough, approximation operator with different prime bases
should provide consistent performance. In other words, the prime
base can not be set as some number which is close to the projection
number. Otherwise, from \eref{DigitExpan} and \eref{VanDerCorput},
we will find that $z_b(m)$ and $m$ may have strong correlation,
e.g. partial linear dependence, that results in a non-uniform sampling
on the sphere. Therefore, we set $b=2$ for rest simulations.

Fig.\ref{fig:DirectPrime} illustrates another exciting phenomenon:
the approximation operator $\cM_{\text{V}_3}^M[\cdot]$ works
better than other two for arbitrarily selected projection number
and prime base. This can be easily understood with following two
considerations: 1) $\cM_{\text{V}_3}^M[\cdot]$ requires higher
order derivatives for smoothness evaluations than other two;
2) many common functions are continuously differentiable with high
order. On the other hand, the performance of
$\cM_{\text{V}_2}^M[\cdot]$ and the one of
$\cM_{\text{M}}^M[\cdot]$ coincide each other perfectly, which
further implies the effectiveness of the proposed method.

\subsection{Decomposition performance}

\begin{figure}[t]
\centering
\includegraphics[width=8.8cm]{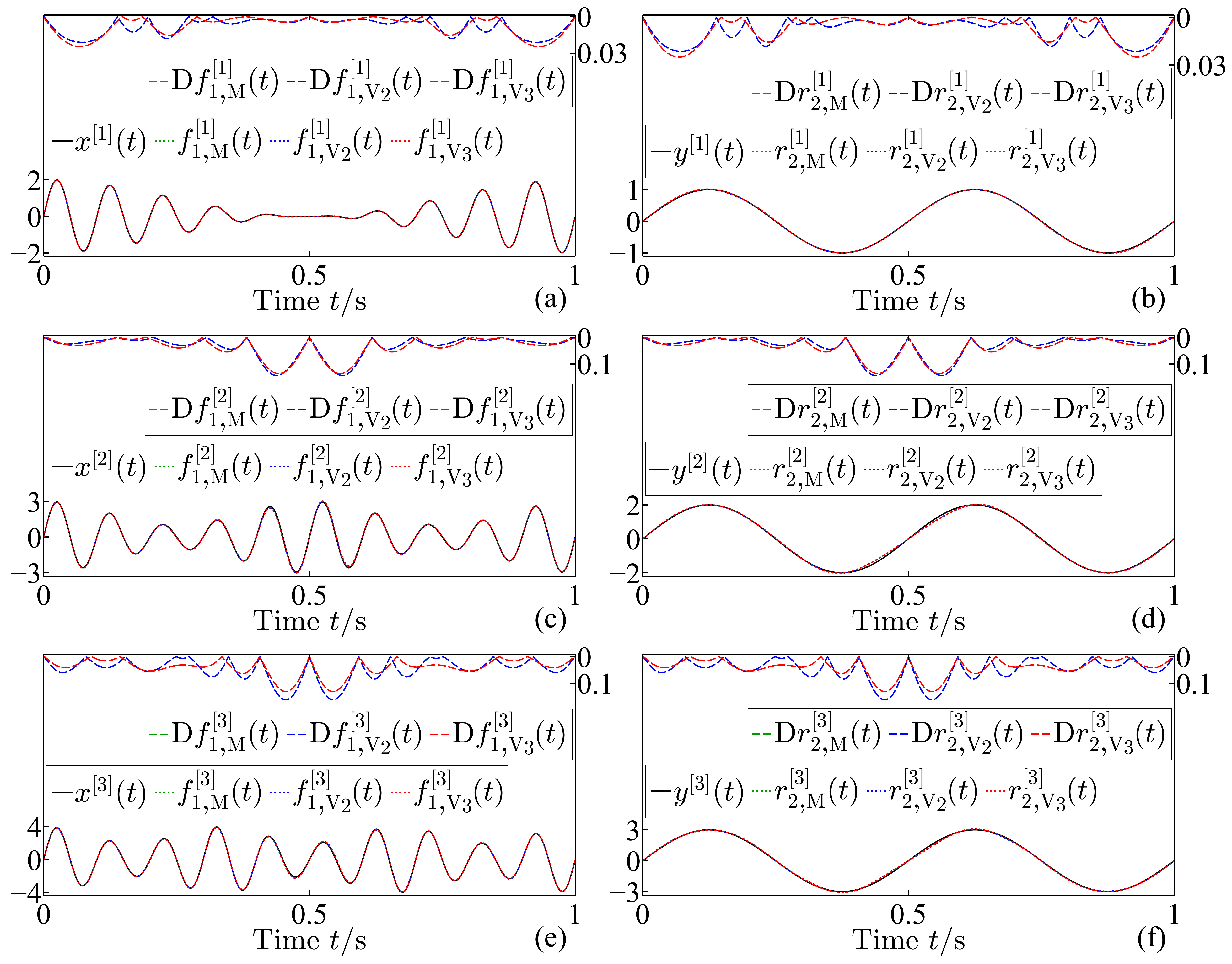}
\caption{\small
Decomposition results of the given VvS $F(t)$ in \eref{Simu1} by
using MEMD and VEMD methods. Left column (a)(c)(e): the original
component $X(t)$, the first decomposed Imfs using MEMD
($F_{1,\text{M}}(t)$), and VEMD with 2nd derivative
($F_{1,\text{V}_2}(t)$) and 3rd derivative ($F_{1,\text{V}_3}(t)$),
and the absolute difference between $X(t)$ and each $F_{1, \cdot}(t)$,
e.g. $\text{D}F_{1,\text{M}}(t) := |X(t) - F_{1,\text{M}}(t)|$;
Right column (b)(d)(f): the original component $Y(t)$, the residual
$R_{2,\text{M}}(t)$, $R_{2,\text{V}_2}(t)$ and $R_{2,\text{V}_3}(t)$,
and the corresponding absolute difference.}
\label{fig:DecompositionRes}
\end{figure}

With the fixed parameters, the proposed VEMD can be applied to
decompose the given VvS $F(t)$ in \eref{Simu1}. The corresponding
decomposition results are shown in Fig.\ref{fig:DecompositionRes}.
As can be seen, all decomposed components by using MEMD or VEMD with
different derivative orders are very similar to the ideal ones. To
numerically distinguish the decomposition performance, the PRDs of
the decomposition results w.r.t each method are listed in
Table-\ref{Tab:PRD}. The numerical results support our observation
again: the VEMD with 3rd derivative order is better than the one
with 2nd derivative order and the MEMD; and the latter two methods
have similar behavior from mathematic point of view.
\begin{table}[t]
\caption{PRD of the decomposition results using MEMD and VEMD}
\centering
\begin{tabular}{cccc}
\\
\hline
       PRD(\%)    &   M  & $\text{V}_2$ & $\text{V}_3$  \\
\hline
 $F_{1,\cdot}(t)$ & 2.95 &    2.95      &     2.65      \\
 $R_{2,\cdot}(t)$ & 3.10 &    3.10      &     2.79      \\
\hline
\end{tabular}
\label{Tab:PRD}
\end{table}

\section{Conclusion}

As a representative data-driven method, EMD can provide finer
time-frequency analysis of any nonlinear and/or nonstationary signal
sampled in (non-)uniform grids. This paper introduced a novel method
to extend the classical EMD for vector-valued signal decomposition.
Different from the existing MEMD, our proposed VEMD obtains the
local mean curve by projecting the envelopes obtained in 1-D projected
space back into the original multi-D space. Since the VEMD does not
require any particular envelope interpolation method, it is more
general and flexible compared to the MEMD for many potential
applications. In addition, with the method introduced in
\cite{Huang2014}, we can generate a time-frequency representation
  of any VvS which can then provide meaningful information for
  further vector-valued data analysis.

\section*{Acknowledgment}
B. Huang is grateful to the Alexander von Humboldt foundation.

\end{document}